\documentclass[11pt]{article}

\usepackage[english]{babel}

\usepackage[center]{caption2}
\usepackage{indentfirst}
\usepackage{amsfonts}
\usepackage{amsmath}
\usepackage{amssymb}
\usepackage{amsbsy, amsthm}

\usepackage[center]{caption2}
\usepackage{amsfonts,amssymb,amsmath,latexsym,amsthm}
\usepackage{multirow}
\usepackage[usenames,dvipsnames]{pstricks}
\usepackage{epsfig}
\usepackage{pst-grad} % For gradients
\usepackage{pst-plot} % For axes
\usepackage[space]{grffile} % For spaces in paths
\usepackage{etoolbox} % For spaces in paths

\newtheorem{dfn}{Definition}

\newtheorem{thm}[dfn]{Theorem}

\newtheorem{lem}[dfn]{Lemma}

\usepackage{latexsym,amsmath,amssymb,amsfonts,epsfig,graphicx,cite,psfrag}
\usepackage{eepic,color,colordvi,amscd}

\def\Bin{\mathrm{Bin}}

\newcommand{\E}{{\mathbb{E}}}
\newcommand{\pr}{{\mathbb{P}}}

\def\QED{\hfill \rule{7pt}{7pt}}

\begin{document}

\title{Partitioning graphs with linear minimum degree}
\author{
Jie Ma\thanks{School of Mathematical Sciences, University of Science and Technology of China, Hefei 230026, China.
Research supported by National Key Research and Development Program of China 2020YFA0713100, National Natural Science Foundation of China grant 12125106, Innovation Program for Quantum Science and Technology 2021ZD0302902, and Anhui Initiative in Quantum Information Technologies grant AHY150200. Email: jiema@ustc.edu.cn.}
\and
Hehui Wu\thanks{
Shanghai Center for Mathematical Sciences,
Fudan University, Shanghai 200438, China.
Research supported by National Natural Science Foundation of China grant 11931006, National Key Research and Development Program of China 2020YFA0713200, and the Shanghai Dawn Scholar Program grant 19SG01. Email: hhwu@fudan.edu.cn.
}
}

\date{}

\maketitle

\begin{abstract}
We prove that there exists an absolute constant $C>0$ such that, for any positive integer $k$, every graph $G$ with minimum degree at least $Ck$ admits a vertex-partition $V(G)=S\cup T$, where both $G[S]$ and $G[T]$ have minimum degree at least $k$, and every vertex in $S$ has at least $k$ neighbors in $T$. This confirms a question posted by K\"uhn and Osthus \cite{KO} and is tight up to a constant factor.
Our proof combines probabilistic methods with structural arguments based on Ore's Theorem on $f$-factors of bipartite graphs.
\end{abstract}

\section{Introduction}
\noindent
There has been extensive research on graph partition problems in graph theory and computer science, due to their various important applications.
%The goal of graph partitioning is often to divide a given graph into several subgraphs or parts, while optimizing some objective functions.
In this paper, we consider the problem of partitioning graphs under minimum degree constraints.

A well known property (often contributed to P. Erd\H{o}s) states that every graph $G$ with minimum degree at least $2k-1$ admits a vertex-partition $V(G)=V_1\cup V_2$ such that
\begin{equation}\label{equ:maxcut}
\mbox{every vertex in $V_i$ has at least $k$ neighbors in $V_{3-i}$ for each $i\in \{1,2\}$.}
\end{equation}
Thomassen \cite{T83} was the first to prove the existence of the least function $f(k)$
such that every graph $G$ with minimum degree at least $f(k)$ admits
a vertex-partition $V(G)=V_1\cup V_2$ satisfying that
\begin{equation}\label{equ:Thom}
\mbox{$G[V_1]$ and $G[V_2]$ both have minimum degree at least $k$.}
\end{equation}
The estimation on $f(k)$ was soon improved by Hajnal in \cite{Haj}.
Finally, Stiebitz \cite{S96} solved a conjecture of Thomassen \cite{T88} by determining the function $f(k)=2k+1$; it is tight as showing by the cliques.
We point out that the analogous partitioning problem for maximum degree was settled by Lov\'asz \cite{Lov} earlier,
while the one for average degree was treated only recently in \cite{CLNWY,WW}.

It is natural to ask if for a graph $G$ with sufficiently large minimum degree, there exists a partition $V(G)=V_1\cup V_2$ such that both \eqref{equ:maxcut} and \eqref{equ:Thom} hold (i.e., every vertex in $V(G)$ has at least $k$ neighbors in each $V_i$).
It turns out that this is impossible even for $k=1$ as shown by the following example of K\"uhn and Osthus \cite{KO}:
Let $\ell$ be any integer and $H=(X,Y)$ be the bipartite graph with $|X|=n$ and $Y=\binom{X}{\ell}$,\footnote{The set $\binom{X}{\ell}$ denotes the family of all subsets of size $\ell$ in $X$.} where $i\in X$ is adjacent to $A\in Y$ if and only of $i\in A$.
It is clear that $H$ has minimum degree $\ell$ and, when provided $n\geq 2\ell-1$, given any bipartition $V(H)=V_1\cup V_2$, one can always find a vertex with none of its neighbors in some $V_i$.\footnote{To see this, note that as $n\geq 2\ell-1$, there exists some $V_j$ with $|V_j\cap X|\geq \ell$; then every vertex $A\in \binom{V_j\cap X}{\ell}$ has all of its $k$ neighbors in $V_j$ and thus zero neighbors in the other partite set.}
On the other hand, K\"uhn and Osthus \cite{KO} proved the following strengthening, which shows that \eqref{equ:Thom} and one side of \eqref{equ:maxcut} can hold simultaneously.

\begin{thm}[K\"uhn-Osthus; Theorem~1 of \cite{KO}]\label{Thm:KO}
For any positive integer $k$, there exists a function $f(k)$ such that any graph $G$ of minimum degree at least $f(k)$ admits a partition $V(G)=S\cup T$, where both $G[S]$ and $G[T]$ have minimum degree at least $k$ and every vertex in $S$ has at least $k$ neighbors in $T$.
\end{thm}

This beautiful result was used in \cite{KO} to derive the analogous theorem for connectivity and then the existence of non-separating structures in
highly-connected graphs. 
It is also related to many major conjectures, i.e., Lov\'asz removable path conjecture \cite{L75}; see discussions in \cite{KO}.
The proof of Theorem~\ref{Thm:KO} is structural, which results in a quadratic bound $f(k)=O(k^2)$.
K\"uhn and Osthus \cite{KO} asked whether this quadratic bound can be replaced by a linear bound on $k$.

In the present paper, we provide an affirmative answer to the above question of K\"uhn and Osthus.
The following is our main result.

\begin{thm}\label{Thm:main}
There exists a constant $C>0$ such that the following holds. For any positive integer $k$, any graph $G$ of minimum degree at least $Ck$ admits a vertex-partition $V(G)=S\cup T$ such that
both $G[S]$ and $G[T]$ have minimum degree at least $k$ and every vertex in $S$ has at least $k$ neighbors in $T$.
\end{thm}

Our proof approach is distinct.
A basic idea is to utilize probabilistic arguments through a random partition.
However, an apparent challenge we face is the need for operations on certain vertices with significant deviations, which can cause the random partition to become fragile and break down.
To overcome this challenge, we derive a useful local structure by applying Ore's Theorem on $f$-factors of bipartite graphs.
Loosely speaking, this theorem yields two vertex-subsets, $X$ and $Y$, that form a ``bipartite graph'' where each vertex in each part has nearly the same number of neighbors in the other part.
By focusing on a local structure similar to the above example given by K\"uhn and Osthus,
we define an appropriate partition that is uniquely determined by a random progress.
Finally, we prove that this partition satisfies our requirements (for more details, see the proof of Theorem~\ref{Thm:dominating}).

Let $G$ be a graph. For subsets $X, Y\subseteq V(G)$ (not necessarily disjoint), let $E_G(X,Y)=\{xy\in E(G):x\in X \mbox{ and } y\in Y\}$ and let $e_G(X,Y)=|E_G(X,Y)|$.
If $X$ consists of a single vertex $x$, then we write as $E_G(x,Y)$ and $e_G(x,Y)$ respectively.
We often drop the subscript when there is no ambiguity from the context.
We define $[n]$ to be the set $\{1,2,...,n\}$ for positive integers $n$.
Throughout the paper, for simplicity we do not try to optimize the constants used in the calculations.

\section{The proof}
\noindent In this section we present the full proof of Theorem~\ref{Thm:main}.
It consists of two parts: the reduction to an alternative statement--Theorem~\ref{Thm:dominating} and the proof of Theorem~\ref{Thm:dominating}.

\subsection{Reduction to Theorem~\ref{Thm:dominating}}
\noindent
In this subsection, we reduce the proof of Theorem~\ref{Thm:main} to the following result.
For a graph $G$, we say a subset $A\subseteq V(G)$ is {\it $k$-dominating} in $G$ if every vertex of $G$ has at least $k$ neighbors in $A$.

\begin{thm}\label{Thm:dominating}
There exists an integer $k_0>0$ such that the following holds.
For any $k\geq k_0$, any graph $G$ of minimum degree at least $50k$ admits a vertex-partition $V(G)=A\cup B$ such that
$A$ is $k$-dominating in $G$ and $G[B]$ has average degree at least $2k$.
\end{thm}

{\noindent \bf Proof of Theorem~\ref{Thm:main} (assuming Theorem~\ref{Thm:dominating}).}
First consider $k\geq k_0$. We claim that in this case, the constant $C$ can be taken to be $C=50$.
Let $G$ be any graph with minimum degree at least $50k$.
By Theorem~\ref{Thm:dominating}, there exists a partition $V(G)=A\cup B$ such that
$A$ is $k$-dominating in $G$ and $G[B]$ has average degree at least $2k$.
There exists a subset $S\subseteq B$ such that the minimum degree of $G[S]$ is at least $k$.
Let $T=A\cup (B\setminus S)$.
Since $A\subseteq T$, it is clear that $T$ is also $k$-dominating in $G$.
Now $V(G)=S\cup T$ is a desired partition of Theorem~\ref{Thm:main}.
To extend this case to all positive integers $k$, it suffices to take $C=50k_0$.
This proves Theorem~\ref{Thm:main}. \QED

\bigskip

Before we give the proof of Theorem~\ref{Thm:dominating}, we state some preliminary tools as follows.
The first one is the classic Chernoff bound for Binomial Distribution.

\begin{lem}[see \cite{AS}]\label{lem:chernoff}
Let $X\sim \Bin(n,p)$ and let $\mu=\E[X]$. For any $\delta\in (0,1)$, we have $\pr(X\geq (1+\delta)\mu)\leq \exp(-\delta^2\mu/3)$ and
$\pr(X\leq (1-\delta)\mu)\leq \exp(-\delta^2\mu/2)$, where $\exp(x)=e^x$.
\end{lem}

We say a family $\mathcal{A}$ of subsets of $[n]$ is {\it monotonically increasing} if $A\in \mathcal{A}$ and $A\subseteq A' \Rightarrow A'\in \mathcal{A}$.
Fix a real $p\in (0,1)$ and consider the probability distribution obtained by choosing each $i\in [n]$ independently with probability $p$.
The following correlation inequality is helpful for probabilistic estimations in the proof.

\begin{lem}[Kleitman's Lemma; see Theorem 6.3.2 of \cite{AS}]\label{lem:Kle}
Let $\mathcal{A}$ and $\mathcal{B}$ be two monotonically increasing families of subsets of $[n]$.
Then we have $\pr(\mathcal{A}\cap \mathcal{B})\geq \pr(\mathcal{A})\cdot \pr(\mathcal{B}).$
\end{lem}

We also need a theorem of Ore \cite{Ore56,Ore57}, which is a generalization of Hall's Theorem \cite{Hall}
as well as a special case of Tutte's $f$-factor theorem \cite{T52} on bipartite graphs.

\begin{thm}[Ore's Theorem \cite{Ore56,Ore57}]\label{Thm:Ore}
Given a bipartite graph $G[V_1, V_2]$, let $f: V(G)\mapsto \mathbb{N}_{\geq 0}$.
Then there is a subgraph $H$ of $G$ such that $d_H(x)= f(x)$ for all $x\in V_1$ and $d_H(y)\le f(y)$ for all $y\in V_2$ if and only if for all $X\subseteq V_1$, we have
\begin{equation}\label{equ:Ore1}
\sum_{x\in X}f(x)\le \sum_{y\in V_2}\min\{f(y),e_G(y,X)\}.
\end{equation}
\end{thm}

\subsection{Proof of Theorem~\ref{Thm:dominating}}
\noindent We now prove Theorem~\ref{Thm:dominating}. Let $k_0$ be sufficiently large and $G$ be a graph of minimum degree $\delta(G)\geq 50k$ where $k\geq k_0$.
We aim to show that there exists a partition $V(G)=A\cup B$ such that
$A$ is $k$-dominating in $G$ and $G[B]$ has average degree at least $2k$.

We first establish a useful structure in the following claim.

\medskip

{\noindent \bf Claim 1.} There exist an integer $t\geq 5k$, non-empty subsets $X, Y\subseteq V(G)$, and a directed spanning subgraph $D$ of $G$
such that the following hold:\footnote{Here $X$ and $Y$ are not necessarily disjoint and it is allowed to have both arcs $u\to v$ and $v\to u$ in $D$.}
\begin{itemize}
\item [(A).] every vertex $v\in V(G)$ has $d_D^+(v)=5k$ and $d_D^-(v)\leq t$,
\item [(B).] every vertex $x\in X$ satisfies that $N_G(x)\setminus N_D^+(x)\subseteq Y$, and
\item [(C).] every vertex $y\in Y$ has $N_D^-(y)\subseteq X$ and $d_D^-(y)\in \{t-1,t\}$.
\end{itemize}

{\noindent \bf Proof of Claim 1.}
We prove this using Theorem~\ref{Thm:Ore}.
Among all oriented spanning subgraphs $D$ of $G$ satisfying that $d_D^+(v)=5k$ for every $v\in V(G)$,
we choose $D$ so that its maximum in-degree $t:=\Delta^+(D)$ is minimized.
An equivalent way to define $t$ is to use the following auxiliary bipartite graph $G^*[V_1, V_2]$,
where $V_1, V_2$ are two disjoint copies of $V(G)$ and $x\in V_1$ is adjacent to $y\in V_2$ in $G^*$ if and only if $xy\in E(G)$.
Define $f: V(G^*)\mapsto \mathbb{N}_{\geq 0}$ be such that $f(x)=5k$ for every $x\in V_1$ and $f(y)=t$ for every $y\in V_2$.
Then $t$ is the minimum integer such that there exists a subgraph $H$ of $G^*$ under the conditions that
\begin{equation}\label{equ:Ore2}
\mbox{$d_H(x)=f(x)$ for every $x\in V_1$ and $d_H(y)\leq f(y)$ for every $y\in V_2$.}
\end{equation}
It is clear that there is a one-to-one correspondence between this subgraph $H$ of $G^*$ and an oriented spanning subgraphs $D(H)$ of $G$
satisfying that for each $v\in V(G)$, $d_{D(H)}^+(v)=d_H(v_1)$ and $d_{D(H)}^-(v)=d_H(v_2)$, where $v_i$ is the corresponding copy of $v$ in $V_i$ for $i=1,2$.

%In particular, there is no such subgraph $H$ if we replace $f(y)=t$ by $f(y)=t-1$ for all $y\in V_2$.
Let us consider the following evolution of \eqref{equ:Ore2}, where the initial setting is that $f(x)=5k$ for every $x\in V_1$ and $f(y)=t-1$ for all $y\in V_2$
and in each following round, we pick one vertex $y\in V_2$ at a time and then change $f(y)=t-1$ to $f(y)=t$.
This terminates when $f(y)=t$ for all $y\in V_2$.
By the minimality of $t$, there does not exist a subgraph $H$ of $G^*$ satisfying \eqref{equ:Ore2} at the beginning of the evolution,
while there does exist such a subgraph $H$ at the end of the evolution.
Now consider the first moment that there exists a subgraph $H$ of $G^*$ satisfying \eqref{equ:Ore2}.
Let $D=D(H)$ be the oriented spanning subgraph of $G$.
By Theorem~\ref{Thm:Ore}, at this moment there must be some non-empty subset $X\subseteq V_1$ which just turned the violation of \eqref{equ:Ore1} to satisfaction.
Since we only picked one vertex (say $y'\in V_2$) in the previous round and enlarged the value of $f(y')$ by one,
this means that at this very moment, \eqref{equ:Ore1} becomes an equation for $X$, namely,
\begin{equation}\label{equ:Ore3}
\sum_{x\in X}f(x)=\sum_{y\in V_2}\min\{f(y),e_G(y,X)\}.
\end{equation}
%Let $D$ be the oriented spanning subgraph of $G$ corresponding to the present subgraph $H$ of $G^*$.
We define
$$Y=\{y\in V_2: |N_D^-(y)\cap X|=f(y)\}.$$
We point out that $Y$ is non-empty.
Indeed, one can easily infer that every $y\in V_2$ satisfies $|N_D^-(y)\cap X|=\min\{f(y),e_G(y,X)\}$;
moreover, the vertex $y'\in V_2$ just turned the value of $\min\{f(y'),e_G(y',X)\}$ from $t-1$ to $t$,
implying that $|N_D^-(y')\cap X|=f(y')$ and thus $y'\in Y$.
It remains to verify that the claim holds for $t, X, Y$ and $D$.
It is easy to see that Item (A) holds by \eqref{equ:Ore2}, and Item (C) follows from the fact that $f(y)\in \{t-1,t\}$ for all $y\in V_2$.
Observe that if $y\in V_2\setminus Y$, then $|N_D^-(y)\cap X|=\min\{f(y),e_G(y,X)\}<f(y)$ and thus $|N_D^-(y)\cap X|=e_G(y,X)$.
This also shows that all edges in $E_G(X,V_2\setminus Y)$ are arcs from $X$ to $V_2\setminus Y$ in $D$,
which implies Item (B). This proves Claim~1.\QED

\bigskip

Write $V=V(G)$ and let $p$ be a real in $(0,1/2]$. We randomly generate a vertex subset $W_p\subseteq X\cup Y$, where each $v\in X\cup Y$ is selected with probability $p$ independent of other vertices.
Define $$L=\{v\in V(G): |N_G(v)\cap (V\setminus W_p)|<k\} \mbox{ ~ and ~ } N_D^+(L)=\bigcup_{v\in L} N_D^+(v).$$
It is evident that
$$\mbox{$(V\setminus W_p)\cup N_D^+(L)$ is a $k$-dominating set in $G$.}$$
Let $c=50$ and for $X\sim \Bin(n,p)$, we write $b_{n,p}(\ell):=\pr(X\geq \ell)$.
If we randomly generate a subset $U_p$ by selecting each vertex $v\in V$ with probability $1-p$ independent of other vertices,
then for each $v\in V$, we have $|N_G(v)\cap U_p|\sim \Bin\big(d_G(v),1-p\big)$;
using coupling and the fact $\delta(G)\geq ck$, we see that
\begin{equation}\label{equ:pr(L)}
\pr(v\in L) \leq \pr(\mbox{$v$ has less than $k$ neighbors in $U_p$})\leq b_{ck,p}((c-1)k).
\end{equation}

We divide the coming proof in two cases according to the value of $t$ given by Claim~1.

\bigskip

{\noindent \bf Case I.} Suppose that $t\leq \frac{1}{50}\exp((c-2)^2k/6c)$, where $c=50$ and $k\geq k_0$ is large.

\medskip

In this case, we choose $p=1/2$.
Then by the equation \eqref{equ:pr(L)} and Lemma~\ref{lem:chernoff}, for each $v\in V$ we have
$$\pr(v\in L)\leq b_{ck,1/2}\big((c-1)k\big)\leq \exp\left(-\frac{(c-2)^2k}{6c}\right).$$
By Claim~1, each vertex $v\in V$ has at most $t$ in-neighbors in $D$, implying that
$$\pr(v\in N_D^+(L))\leq \sum_{u\in N_D^-(v)} \pr(u\in L)\leq t\cdot \exp\left(-\frac{(c-2)^2k}{6c}\right)\leq \frac{1}{50}.$$
Let $A:=(V\setminus W_{1/2})\cup N_D^+(L)$ and $B:=V\setminus A=W_{1/2}\setminus N_D^+(L).$
We have seen that $A$ is $k$-dominating in $G$.
It remains to consider the average degree of $G[B]$.

For each $v\in X\cup Y$ we have $\pr(v\in B)\leq \frac12$, so $\E[|B|]\leq \frac{|X\cup Y|}{2}.$
By Items (B) and (C) of the claim, each vertex in $X\cup Y$ has at least $\min\{t-1, \delta(G)-k\}\geq 5k-1$ neighbors in $X\cup Y$.
Using the above estimations, we can derive that
\begin{align*}
\E[e(G[B])]&=\sum_{uv\in E(G[X\cup Y])} \pr(\{u,v\}\subseteq B)\geq \sum_{uv\in E(G[X\cup Y])} \bigg(\pr(\{u,v\}\subseteq W_{1/2})-\pr(\{u,v\}\cap N_D^+(L)\neq \emptyset)\bigg)\\
&\geq \sum_{uv\in E(G[X\cup Y])} \bigg(\pr(\{u,v\}\subseteq W_{1/2})-\left(\pr(u\in N_D^+(L))+\pr(v\in N_D^+(L))\right)\bigg)\\
&\geq \sum_{uv\in E(G[X\cup Y])}\left(\frac14-\frac{1}{25}\right)\geq \frac{21}{100}\cdot (5k-1)\cdot \frac{|X\cup Y|}{2} \geq k \cdot \E[|B|].
\end{align*}
That is, $\E[e(G[B])-k|B|]\geq 0$.
Therefore with positive probability, there exists a desired vertex-partition $V(G)=A\cup B$.
We have finished the proof of Case~I when $t\leq \frac{1}{50}\exp\left(\frac{(c-2)^2k}{6c}\right)$.

\bigskip

{\noindent \bf Case II.} Suppose that $t\geq \frac{1}{50}\exp((c-2)^2k/6c)$, where $c=50$ and $k\geq k_0$ is large.

\medskip

We need to choose a suitable $p\in (0,1/2]$ and consider a modification of the vertex-partition given by the previous case (to be more precise, we will replace $B$ with a proper subset).
Let $$S=\{v\in X\cap W_p: |N_G(v)\cap Y\cap W_p|\geq 5k\}.$$
In the rest of the proof, we define
$$ B:=\big((Y\cap W_p)\cup S\big)\setminus N_D^+(L) \mbox{ ~ and ~ } A:=V\setminus B.$$
Since $B\subseteq W_p\setminus N_D^+(L)$, we see that $A\supseteq (V\setminus W_p)\cup N_D^+(L)$ is $k$-dominating in $G$.
Again, we aim to show that the expected average degree of $G[B]$ is at least $2k$, which would imply that with positive probability,
the partition $V(G)=A\cup B$ are as desired, thus completing the proof.

Let us point out that the probability space $\{W_p\}$ considered here is the collection of all subsets of $X\cup Y$,
so the event $v\in W_p$ (for any $v\in X\cup Y$) and the event $x\in S$ (for any $x\in X$) both are monotonically increasing.
Then using Lemma~\ref{lem:Kle}, for any $v\in X\cup Y$ and $x\in X$ we have
\begin{equation}\label{equ:Kle}
\pr(v\in W_p|x\in S)\geq \pr(v\in W_p)=p.
\end{equation}

We first choose a suitable probability $p\in (0,1/2)$ in the following claim.

\medskip

{\noindent \bf Claim 2.} There exists a real $p\in (0,1/2)$ such that $$p\cdot b_{(c-5)k,p}(5k)=\frac{5k}{t} \mbox{ ~ and ~ } b_{ck,p}((c-1)k)\leq \frac{1}{t^3}.$$

{\noindent \bf Proof of Claim 2.}
Let $g(p)=p\cdot b_{(c-5)k,p}(5k)$ be a function with variable $p$.
It is easy to see that $g(p)$ is an increasing continuous function with $g(0)=0$ and $g(1/2)=1/2\cdot b_{(c-5)k,1/2}(5k)\geq 1/4$.
Since $0<\frac{5k}{t}\leq 250k/\exp((c-2)^2k/6c)<1/4$,
there exists a unique $p\in (0,1/2)$ satisfying $g(p)=\frac{5k}{t}$.

Note that $t\geq \frac{1}{50}\exp((c-2)^2k/6c)$, so $e^k\leq (50t)^{6c/(c-2)^2}$.
Using basic properties of the binomial distribution,
we see $\frac{5k}{t}=p\cdot b_{(c-5)k,p}(5k)\geq p^{5k+1}$, which implies that $p\leq (5k/t)^{1/(5k+1)}$,
and $$b_{ck,p}((c-1)k)\leq \binom{ck}{(c-1)k} p^{(c-1)k}\leq \left(\frac{ec}{c-1}\right)^{(c-1)k}p^{(c-1)k}\leq e^{ck}p^{(c-1)k},$$
where the last inequality holds because $(1+1/x)^x$ increases to $e$ as $x$ goes to infinity.
Therefore
\begin{align*}
\frac{b_{ck,p}((c-1)k)}{p\cdot b_{(c-5)k,p}(5k)}\leq \frac{e^{ck}p^{(c-1)k}}{p^{5k+1}}\leq (50t)^{\frac{6c^2}{(c-2)^2}} \left(5k/t\right)^{\frac{(c-6)k-1}{5k+1}}\leq \frac{1}{5kt^2},
\end{align*}
where the last inequality holds as $c=50$ and $t\gg k\geq k_0$ is large.
Finally, this shows that $b_{ck,p}((c-1)k)\leq p\cdot b_{(c-5)k,p}(5k)\cdot \frac{1}{5kt^2}=\frac{5k}{t}\cdot \frac{1}{5kt^2}=\frac{1}{t^3},$
proving Claim~2.
\QED

\bigskip

Consider any vertex $x\in X$.
By Claim~1 we have $|N_G(x)\cap Y|\geq d_G(x)-5k\geq (c-5)k$.
Recall the definition of $S$. Then Claim~2 shows that for $x\in X$,
\begin{equation}\label{equ:p(S)}
\pr(x\in S)\geq p\cdot b_{(c-5)k,p}(5k)=\frac{5k}{t}.
\end{equation}
For any vertex $v\in V$, by \eqref{equ:pr(L)} and Claim~2 we have $\pr(v\in L)\leq b_{ck,p}((c-1)k)\leq \frac{1}{t^3}$.
Recall that the maximum in-degree of $D$ is at most $t$ (i.e., Item (A) of Claim~1).
So we can obtain
\begin{equation}\label{equ:p(N)}
\pr(v\in N_D^+(L))\leq \sum_{u\in N_D^-(v)} \pr(u\in L)\leq t\cdot \frac{1}{t^3}=\frac{1}{t^2}.
\end{equation}
%This tells that the probability of any vertex being in $N_D^+(L)$ is relatively small.
%So the subset $B=\big((Y\cap W_p)\cup S\big)\setminus N_D^+(L)$ should behave just like $(Y\cap W_p)\cup S$, the latter of which has high average degree.
%To prove this rigorously,
Then, for any $x\in X$ and $y\in Y$ with $xy\in E(G)$, we have
\begin{align*}\label{equ:error}
\pr(\{x,y\}\subseteq B)&\geq \pr(x\in S \wedge y\in W_p)-\pr(\{x,y\}\cap N^+_D(L)\neq \emptyset)\\
&\geq \pr(x\in S)\cdot \pr(y\in W_p|x\in S)-\left(\pr(x\in N_D^+(L))+\pr(y\in N_D^+(L))\right)\\
&= \pr(x\in S)\cdot \pr(y\in W_p|x\in S)\cdot\left(1-\frac{\pr(x\in N_D^+(L))+\pr(y\in N_D^+(L))}{\pr(x\in S)\cdot \pr(y\in W_p|x\in S)}\right)\\
&\geq \pr(x\in S)\cdot \pr(y\in W_p|x\in S)\cdot \left(1-\frac{2/t^2}{(5k/t)\cdot p}\right)\\
&\geq 0.99\cdot \pr(x\in S)\cdot \pr(y\in W_p|x\in S),
\end{align*}
where the second last inequality follows by \eqref{equ:Kle}, \eqref{equ:p(S)} and \eqref{equ:p(N)},
and the last inequality holds because $p\geq p\cdot b_{(c-5)k,p}(5k)=\frac{5k}{t}$ (using Claim~2).
Using this inequality, we can derive the following (note that $X$ and $Y$ may overlap)
\begin{align*}
4\cdot \E[e(G[B])]\geq &\left(\sum_{x\in X}\sum_{y\in Y\cap N_G(x)}+\sum_{y\in Y}\sum_{x\in X\cap N_G(y)}\right) \pr(\{x,y\}\subseteq B)\\
\geq &\left(\sum_{x\in X}\sum_{y\in Y\cap N_G(x)}+\sum_{y\in Y}\sum_{x\in X\cap N_G(y)}\right)0.99\cdot \pr(x\in S)\cdot \pr(y\in W_p|x\in S)\\
\geq &\sum_{x\in X}0.99\cdot \pr(x\in S)\cdot \bigg(\sum_{y\in Y\cap N_G(x)}\pr(y\in W_p|x\in S)\bigg)+ \sum_{y\in Y} \sum_{x\in X\cap N_G(y)}0.99\cdot \frac{5k}{t}\cdot p\\
\geq &\sum_{x\in X}0.99\cdot \pr(x\in S)\cdot 5k+ \sum_{y\in Y} (t-1)\cdot \frac{0.99\cdot 5k}{t}\cdot p\\
\geq & 4k\cdot \E[|S|]+4k\cdot \E[|Y\cap W_p|]\geq 4k\cdot \E[|B|],
\end{align*}
where the third inequality follows from \eqref{equ:Kle} and \eqref{equ:p(S)},
the fourth inequality holds by the definition of $S$ and
the fact from Item~(C) of Claim~1 that every $y\in Y$ has at least $t-1$ neighbors in $X$,
and the last inequality holds by the definition of $B$.
This proves that $\E[e(G[B])-k|B|]\geq 0$, completing the proof of Case~II (and thus the proof of Theorem~\ref{Thm:dominating}).
\QED

\bigskip

\bigskip

{\noindent \bf Acknowledgment.} The authors would like to thank Shengjie Xie for very helpful discussions and for his careful reading on a preliminary draft.

\end{document}